\documentclass[12pt,a4paper,reqno]{amsart}
\usepackage{fullpage}
\usepackage{amsmath,amssymb}
\usepackage{braket}
\usepackage{verbatim}
\usepackage{lscape}
\usepackage{newtxtext,newtxmath}
\usepackage{hyperref}
\usepackage{url}
\usepackage{xcolor}

\newcommand{\Odi}[1]{O\left(#1\right)}
\newcommand{\dx}{\,\mathrm{d}}
\newcommand{\e}{\mathrm{e}}
\newcommand{\ii}{\mathrm{i}}
\newcommand{\m}{\mathfrak{m}}

\newcommand{\C}{\mathcal{C}}

\newcommand{\E}{\mathbb{E}}
\newcommand{\N}{\mathbb{N}}
\newcommand{\R}{\mathbb{R}}

\newcommand{\Sums}{\mathfrak{S}}

\DeclareMathOperator\supp{supp}

\newtheorem{Theorem}{Theorem}[section]
\newtheorem{Proposition}{Proposition}[section]
\newtheorem{Corollary}{Corollary}[section]
\newtheorem{Lemma}{Lemma}[section]
\newtheorem{Remark}{Remark}[section]

\allowdisplaybreaks

\author{Alessandro Gambini, Remis Tonon,
        and Alessandro Zaccagnini  \\ with an addendum by Jacques Benatar and Alon Nishry}

\title{Signed harmonic sums of integers \\ with $k$ distinct prime factors}

\date{\today}

\begin{document}

\begin{abstract}
We give some theoretical and computational results on ``random''
harmonic sums with prime numbers, and more generally, for integers
with a fixed number of prime factors.
\end{abstract}

\maketitle

\small


\normalsize
\section{Introduction and general setting}

It is well known that the harmonic series restricted to prime numbers
diverges, as the harmonic series itself.
This was first proved by Leonhard Euler in 1737 \cite{Euler1737}, and
it is considered as a landmark in number theory.
The proof relies on the fact that
\[
  \sum_{n = 1}^N \frac1n
  =
  \log N
  +
  \gamma
  +
  \Odi{1/N},
\]
where $\gamma \simeq 0.577215\dots$ is the Euler--Mascheroni constant.
The corresponding result for primes is one of the formulae proved
by Mertens, namely
\[
  \sum_{p\le N}\frac1p=\log\log N+A+\Odi{\frac{1}{\log N}},
\]
where $A\simeq 0.2614972\ldots$ is the Meissel--Mertens constant.
It is also referred to as Hadamard--de la Vall\'ee-Poussin constant
that appears in Mertens' second theorem.

Recently, Bettin, Molteni and Sanna \cite{BettinMS2018b} studied the
random harmonic series
\begin{equation}
\label{rand-harm-sum}
  X := \sum_{n=1}^{\infty}\frac{s_n}{n},
\end{equation}
where $s_1$, $s_2,\ldots$ are independent uniformly distributed random
variables in $\{-1,+1\}$.
Based on the previous work by Morrison \cite{Morrison1995,Morrison1998}
and Schmuland \cite{Schmuland2003}, they proved the almost sure
convergence of \eqref{rand-harm-sum} to a density function $g$,
getting lower and upper bounds of the minimum of the distance of a
number $\tau\in\mathbb{R}$ to a partial sum $\sum_{n=1}^{N}s_n/n$.
In 1976 Worley studied the same problem in terms of upper bound of \eqref{rand-harm-sum} both in the case $\tau=0$ (see \cite{Worley1976a}) and for a generic $\tau\in\mathbb{R}$ (see \cite{Worley1976b}); his approach is not probabilistic but he has achieved an upper bound comparable to that of \cite{BettinMS2018b}.
For further references, see also Bleicher and {Erd\H os}
\cite{BleicherE1975,BleicherE1976}, where the authors treated the
number of distinct subsums of $\sum_1^N 1/n$, which corresponds to
taking $s_i$ independent uniformly distributed random variables in
$\{0,1\}$.
A more complete list of references can be found in
\cite{BettinMS2018b}.

The purpose of this paper is firstly to show that basically the same
results hold for a general sequence of integers under some suitable,
and not too restrictive, conditions; moreover, that a stronger result
can be reached if we restrict to integers with exactly $k$ distinct
prime factors.

Although Bettin, Molteni and Sanna \cite{BettinMS2018b} treat both the lower
bound and the upper bound, we are mainly interested in the upper bound
using a probabilistic approach.
As we will see, in the cases that we treat, we will not be able to say
anything about the lower bound, except in terms of numerical
computations.

We will use a consistent notation with the previous works by Bettin,
Molteni and Sanna \cite{BettinMS2018a}, \cite{BettinMS2018b}, Crandall
\cite{Crandall2008} and Schmuland \cite{Schmuland2003}.

\subsection{General setting of the problem}

We denote by $\N$ the set of positive integers.
Let $(a_n)_{n \in \N}$ be a strictly decreasing sequence of positive
real numbers such that
\begin{equation}
\label{general-seq}
  \lim_{n \to +\infty} a_n
  =
  0
  \qquad\text{and}\qquad
  \sum_{n \ge 1} a_n
  =
  + \infty.
\end{equation}
Notice that
\[
  \sum_{n \ge 1} (-1)^n a_n
\]
converges (not absolutely) by Leibniz's rule.
Hence, by Riemann's theorem, given
$\lambda$, $\Lambda \in [-\infty, +\infty]$ with
$\lambda \le \Lambda$, we can arrange the choice of the signs
$s_n = s_n(\lambda, \Lambda) \in \{ -1, 1 \}$, in such a way that
\[
  \liminf_{N \to +\infty}
    \sum_{n \le N} s_n a_n
  =
  \lambda
  \qquad\text{and}\qquad
  \limsup_{N \to +\infty}
    \sum_{n \le N} s_n a_n
  =
  \Lambda.
\]
As we said above, we are mainly interested in prime numbers, so we
introduce some further reasonable hypotheses on the sequence $a_n$:
we assume that $b_n = a_n^{-1} \in \N$, so that $b_n$ is strictly
increasing, and that
\begin{equation}
\label{hyp-bn}
  n \le b_n \le n B(n),
\end{equation}
where $B(n)=n^{\beta(n)}$, with $\beta$ a real-valued decreasing function such that $\beta(n)=o(1)$.
In order to prove Proposition \ref{prop32} below, we will assume a
more restrictive condition on $\beta$, that is
\begin{equation}
\label{eq:conditionbeta}
\beta(n)\le \frac {1}{8\log\log n}
\qquad\text{for sufficiently large $n$}.
\end{equation}
Actually, this assumption is not strictly necessary and we will
discuss this in Remark \ref{rk:beta}. Nevertheless, since the series
$\sum a_n$ must diverge, this condition is not too restrictive, and
besides it is satisfied by most of the interesting sequences, like
arithmetic progressions, the one of primes, and primes in arithmetic
progressions.

Let us introduce some more notation: we consider the set
\begin{equation}
\label{def-sums}
  \Sums(N)
  =
  \Bigl\{
    \sum_{n \le N} s_n a_n
    \colon
    s_n \in \{ \pm  1 \}
    \text{ for }
    n \in \{1,  \dots, N \}
  \Bigr\},
\end{equation}
and, for a given $\tau \in \R$, we set
\[
  \m_N(\tau)
  =
  \min \bigl\{ \vert S_N - \tau \vert \colon S_N \in \Sums(N) \bigr\}.
\]
In other words, for a given $N \in \N$, the goal is to find the choice
of signs such that $\vert S_N -\tau \vert$ attains its minimum value.
Finally, we define the random variable
\[
  X_N := \sum_{n=1}^{N} s_n a_n,
\]
where the signs $s_n$ are taken uniformly and independently at random
in $\{-1,1\}$.
We will study its small scale distribution.
With a slight abuse of notation, we denote by $s_n$ both the signs
in the definition \eqref{def-sums} and the random variables in the
definition above.

\subsection{Results}

For ease of comparison with the results in Bettin, Molteni and Sanna
\cite{BettinMS2018b}, we now state our main results in the following
form, even though more precise versions of them are to be found within
the paper.

\begin{Theorem}\label{genthm}
Let $\beta$ satisfy \eqref{eq:conditionbeta}.
Then there exists $C>0$ such that for every $\tau\in\R$ we have
\[
\m_N(\tau)<\exp(-C\log^2 N)
\]
for all sufficiently large $N$ depending on $\tau$.
\end{Theorem}

\begin{Theorem}\label{primetheorem}
Let $(b_n)_{n\in\N}$ be the sequence of integers having exactly
$k$ distinct prime factors. Then, for every $\tau\in\R$ and for all
sufficiently large $N$ depending on $\tau$, we have
\[
\m_N(\tau)<\exp(-f(N)),
\]
where $f$ is any function satisfying
\[
f(N) = o\left(N^{1/(2k+1)-\varepsilon}\right).
\]
\end{Theorem}

\begin{Remark}
We emphasize the fact that the estimate obtained in  Theorem \ref{primetheorem} holds uniformly for every $\tau\in\R$ in any fixed compact set.
\end{Remark}

\begin{Corollary}[J.~Benatar and A.~Nishry]\label{BN}
For any fixed $\tau \in \R$, $\varepsilon > 0$ and any sufficiently
large $N$ there exists a choice of signs
$(s_n)_{n \leq N}\in \left\{-1,1 \right\}^N$, such that
\[
  \left| \sum_{n \leq N} \frac{s_n}{n} - \tau\right|
  \ll_{\tau,\varepsilon}
  \exp \left( -N^{1/3-\varepsilon}\right).
\]
\end{Corollary}

We collect some numerical results for $k = 1$ in Tables~\ref{Fig2},
\ref{Fig3} and \ref{Fig5}. The sequence of  Tables \ref{Fig2} and \ref{Fig3} appears in \href{https://oeis.org/A332399}{OEIS A332399}: see \cite{CGTZ}.

\smallskip
\noindent{\textbf{Acknowledgements.}}
We thank Sandro Bettin and Giuseppe Molteni for many conversations on
the subject, and Mattia Cafferata for his help in computing the tables
at the end of the present paper.
We also warmly thank Jacques Benatar and Alon Nishry for their
fruitful suggestions which improved our paper, for providing us
references and for letting us include their proof of Corollary
\ref{BN} in this paper.
R.~Tonon and A.~Zaccagnini are members of the INdAM group GNSAGA,
which partially funded their participation to the Second Symposium on
Analytic Number Theory in Cetraro, where some of this work was done.

\section{Lemmas}

In this section we study some properties of the general sequence
defined in \eqref{general-seq}, using the classical notation:
$\E[X]$ denotes the expected value of a random variable $X$,
$\mathbb{P}(E)$ the probability of an event $E$.
For each continuous function with compact support $\Phi \in \C_c(\R)$ we
denote by $\widehat{\Phi}$ its Fourier transform defined as follows:
\[
  \widehat{\Phi}(x)
  :=
  \int_{\R} \Phi(y) \, \e^{-2 \pi \ii x y} \dx y.
\]
We are actually interested in smooth functions, because the smoothness
of the density of any random variable $X$ is related to the decay at
infinity of its characteristic function, defined precisely by its
Fourier transform.

For each $N \in \N \cup \{ \infty \}$, for any $x \in \R$ and for any
sequence satisfying \eqref{general-seq}, we also define the product
\[
  \varrho_N(x)
  :=
  \prod_{n=1}^N\cos(\pi xa_n)
  \qquad\text{and}\qquad
  \varrho(x)
  :=
  \varrho_{\infty}(x).
\]

We begin with the following lemma, which is a more general version of
Lemma 2.4 from \cite{BettinMS2018b}.

\begin{Lemma}
\label{media_phi}
We have
\[
  \E[\Phi(X_N)]
  =
  \int_{\R} \widehat{\Phi}(x) \varrho_{N}(2x) \dx x
\]
for all $\Phi \in \C_c^1(\R)$.
\end{Lemma}

\proof
By the definition of expected value we have
\[
  \E[\Phi(X_N)]
  =
  \frac{1}{2^N}
  \sum_{s_1,\ldots,s_N\in\{-1,1\}}
    \Phi \left( \sum_{n=1}^N s_n a_n \right).
\]
Using the inverse Fourier transform we get
\begin{align*}
  \E[\Phi(X_N)]
  &=
  \frac{1}{2^N}
  \sum_{s_1,\ldots,s_N\in\{-1,1\}}
    \int_{\R} \widehat{\Phi}(x)
      \exp\left( 2 \pi \ii x \sum_{n=1}^N s_n a_n \right) \dx x \\
  &=
  \int_{\R} \widehat{\Phi}(x)
    \frac{1}{2^N} \sum_{s_1,\ldots,s_N\in\{-1,1\}}
      \exp\left(2 \pi \ii x \sum_{n=1}^N s_n a_n \right) \dx x.
\end{align*}
Exploiting the fact that
$\e^{\ii \alpha} + \e^{-\ii \alpha} = 2 \cos(\alpha)$,
we have
\[
  \sum_{s_1,\ldots,s_N\in\{-1,1\}}
    \exp\left(2\pi i x\sum_{n=1}^Ns_na_n\right)
  =
  \frac12
  \sum_{s_1,\ldots,s_N\in\{-1,1\}}2
  \cos\left(2\pi x\sum_{n=1}^Ns_na_n\right).
\]
Finally, taking advantage of Werner's trigonometric identities, we obtain
\[
  \E[\Phi(X_N)]
  =
  \int_{\mathbb{R}}\widehat{\Phi}(x)\varrho_N(2x) \dx x.
\qedhere
\]
\endproof

We will need also a generalisation of Lemma 2.5 from \cite{BettinMS2018b}, which is the following

\begin{Lemma}
\label{lem:3.2}
For all $N \in \N$ and $x\in [0,\sqrt{N}]$ we have
\[
  \varrho_N(x)
  =
  \varrho(x) \left(1 + \Odi{x^2/N}\right).
\]
\end{Lemma}

\proof
We recall that $a_n$ is defined as in \eqref{general-seq} and
satisfies \eqref{hyp-bn}.
In particular $a_n = \Odi{1 / n}$, so that the same argument in the
proof of Lemma 2.5 of \cite{BettinMS2018b} holds.
\endproof

Let us now define, for every positive integer $N$ and any real
$\delta$ and $x$ the set
\begin{equation*}
  \mathcal{S} \bigl(N, \delta, x, (a_n)_{n \ge 1} \bigr)
  :=
  \{ n \in \{1,\ldots,N\} \colon \Vert x \, a_n \Vert \ge \delta\},
\end{equation*}
where $\Vert\cdot\Vert$ denotes the distance from the nearest integer.
For brevity, we sometimes drop the dependence on the sequence
$(a_n)_{n \ge 1}$.

\begin{Lemma}
\label{stima_rho_N}
For all $N \in \N$ and for all $x,\delta\ge 0$ we have
\begin{equation*}
  \vert \varrho_N(x) \vert
  \le
  \exp\left(-\frac{\pi^2\delta^2}{2}\cdot\# \mathcal{S}(N,\delta,x)\right).
\end{equation*}

\end{Lemma}
\proof
It is a straightforward consequence of the inequality
\[
  \vert\cos(\pi x)\vert
  \le
  \exp\left(-\frac{\pi^2\Vert x\Vert^2}{2}\right).
\qedhere
\]
\endproof

\begin{Lemma}
\label{lemma34}
For any $N \in \N$, $x \in \R$ and $0<\delta<1/2$ we have
\[
  \frac{N}{2}
  -
  D(N, y(\delta), x)
  <
  \# \mathcal{S}(N,\delta,x)
  <
  N
  -
  D(N, y(\delta)/2, x),
\]
where
\[
D(N, y, x) =
D\left(N,y, x, (b_n)_{n \ge 1}\right)
:= \sum_{x - y < m < x + y}
  \sum_{\substack{b_n \vert m \\ N / 2 \le n \le N}} 1
\]
and $y(\delta) := \delta N B(N)$.
\end{Lemma}

\proof
As in Lemma 3.3 of \cite{BettinMS2018b}, we observe that
\[
  \frac{N}{2} - T(N, \delta, x)
  <
  \#\mathcal{S}(N,\delta,x)
  <
  N - T(N, \delta, x),
\]
where
\[
  T(N, \delta, x)
  :=
  \#\{ n \in \N \cap [N / 2, N] \colon \Vert x a_n \Vert < \delta\}.
\]
Now, recalling that $a_n = 1 / b_n$, we have
\begin{align*}
  T(N, \delta, x)
  &=
  \#\{n\in\mathbb{N}\cap[N/2,N]:\exists\ell\in\mathbb{N},\ \ell-\delta<xa_n<\ell+\delta\}
  \\
  &=
  \#\{n\in\mathbb{N}\cap[N/2,N]:\exists\ell\in\mathbb{N},\ x-\delta b_n<\ell b_n<x+\delta b_n\}.
\end{align*}

From our hypothesis \eqref{hyp-bn} we know that $b_n \le N B(N)$; then
\begin{align*}
  T(N, \delta, x)
  &<
  \#\{n\in\mathbb{N}\cap[N/2,N] \colon
    \exists\ell\in\mathbb{Z},\ x-y(\delta)<\ell b_n<x+y(\delta)\} \\
  &= D(N,y(\delta),x).
\end{align*}
This proves the lower bound; the upper bound follows with the same
argument.
\endproof

\begin{Proposition}
\label{prop32}
Let $A$ be a fixed positive constant and, for $N$ sufficiently large,
\begin{equation*}
\beta(N)\le \frac {1}{8\log\log N}.
\end{equation*}
Then there exists $C'>0$ such that
$\vert\varrho_N(x)\vert<x^{-A}$ for all sufficiently large positive
integers $N$ and for all $x\in[N,\exp(C'(\log N)^2)]$.
\end{Proposition}

\proof
The proof follows along the same lines as Proposition 3.2 of \cite{BettinMS2018b}: we
take
\[
\overline{\delta}=\frac{2\sqrt{2A\log x}}{\pi}N^{-1/2}
\quad \text{and} \quad
x\in\bigg[N,\exp\left(\frac{\pi^2 N}{32A}\right)\bigg),
\]
so that $0<\overline{\delta}<1/2$ and $y(\overline{\delta}) = \overline{\delta} N B(N)<x$.

By Lemmas \ref{stima_rho_N} and \ref{lemma34}, if we show that $D(N,y(\overline{\delta}),x)<N/4$, then we get
$\vert\varrho_N(x)\vert<1/x^A$. Considering that $b_n$ is a sequence of positive integers, we use Rankin's trick with $w\in(1/4,1/2)$ and Ramanujan's result on $\sigma_{-s}(n)$ \cite{Ramanujan1915} to obtain
\begin{align*}
  D(N,y(\overline{\delta}),x)
&<\frac{4}{\pi}\sqrt{2AN\log x}\, B(N)\cdot\max_{m\le 2x}\sum_{\substack{b_n\vert m \\ N/2\le n\le N}}1 \notag \\
&< \frac{4}{\pi}\sqrt{2AN\log x} \, B(N) \cdot \max_{m\le 2x} \sum_{\substack{k \vert m \\ N/2\le k\le N B(N)}} 1 \notag \\
&\le \frac{4}{\pi}\sqrt{2AN\log x} \, B(N) \cdot \max_{m\le 2x} \sum_{\substack{k \vert m \\ N/2\le k\le N B(N)}} \left(\frac{N B(N)}{k}\right)^w \notag \\
&= \frac{4}{\pi}N^{\frac12+w} B(N)^{1+w} \sqrt{2A\log x}\, \cdot\max_{m\le 2x}\sum_{\substack{k \vert m \\ N/2\le k\le N B(N)}} k^{-w} \notag \\
&\le \frac{4}{\pi}N^{\frac12+w} B(N)^{1+w}\sqrt{2A\log x}\,\cdot\max_{m\le 2x} \sigma_{-w}(m) \notag \\
&< \frac{4}{\pi}N^{\frac12+w} B(N)^{1+w} \sqrt{2A\log x}\, \cdot\exp\left(C_1\frac{(\log 2x)^{1-w}}{\log\log 2x}\right), 
\end{align*}
where $C_1$ is the constant of Ramanujan's theorem, as it is stated in Lemma 3.4 of \cite{BettinMS2018b}.

Let $w=w(x):=1/2-\varphi(x)$, where $\varphi$ is a positive decreasing function that we will choose later.
Then we have
\[
B(N)^{1+w} = \exp\left(\left(\frac{3}{2}-\varphi(x)\right)\beta(N)\log N\right),
\]
and so we would be done if we showed that
\[
N^{1-\varphi(x)+(3/2-\varphi(x))\beta(N)} \sqrt{\log x} \cdot\exp\left(C_1\frac{(\log 2x)^{1/2+\varphi(x)}}{\log\log 2x}\right) = o(N),
\]
that is
\begin{equation*}
\sqrt{\log x} \cdot\exp\left(C_1\frac{(\log 2x)^{1/2+\varphi(x)}}{\log\log 2x}\right) = o(N^{\varphi(x)+(\varphi(x)-3/2)\beta(N)}).
\end{equation*}
Hence we must have
\[
\varphi(x)+(\varphi(x)-3/2)\beta(N) > 0,
\]
that is
\begin{equation*}
\beta(N)<\frac{\varphi(x)}{3/2-\varphi(x)} \approx \frac{2}{3}\varphi(x).
\end{equation*}
Since $\varphi$ is decreasing and we want to maintain the same range for $x$ as in \cite{BettinMS2018b}, that is $x\in \left[N,\exp \left(C'(\log N)^2\right)\right]$, we need to have
\[
  \beta(N) \lesssim \frac{2}{3}\varphi\left(\exp \big(C'(\log N)^2\big)\right).
\]
Let us take $\varphi(x) = (\log\log 2x)^{-1}$ and $\beta(N)$ such that for $x\in \left[N,\exp \left(C'(\log N)^2\right)\right]$ it holds
\begin{equation}\label{eq:beta23}
\beta(N) \le \frac{2}{3J}\varphi(x) = \frac{2}{3J} \frac{1}{\log\log 2x},
\end{equation}
where $J\in\R$, $J>1$.
Then we would achieve our goal if we showed that
\[
\sqrt{\log x} \cdot\exp\left(C_1\e\frac{(\log 2x)^{1/2}}{\log\log 2x}\right) = o\left(\exp\left(\left(1-\frac{1}{J}+o(1)\right)\frac{\log N}{\log\log 2x}\right)\right),
\]
that is
\[
\exp\left(C_1\e\frac{(\log 2x)^{1/2}}{\log\log 2x}-\left(1-\frac{1}{J}+o(1)\right)\frac{\log N}{\log\log 2x}+\frac12\log\log x\right) = o(1).
\]
This condition is equivalent to
\[
C_1\e\frac{(\log 2x)^{1/2}}{\log\log 2x}-\left(1-\frac{1}{J}+o(1)\right)\frac{\log N}{\log\log 2x}+\frac12\log\log x \to -\infty.
\]
Taking into account the ranges for $x$, we see that it is sufficient
to have
\[
\frac{1}{\log\log N}\left[C_1 \sqrt{C'}\,\e \log N (1+o(1))-\left(1-\frac{1}{J}\right)\log N+O\left((\log\log N)^2\right) \right]\to -\infty.
\]
We recall that, by our choice of $x$ and $N$, we have $\log\log x
\asymp \log\log N$.
Hence, we just need to take $C'$ sufficiently small, in a way that 
\begin{equation}
\label{eq:C}
C' < \left(\frac{J-1}{C_1\e J}\right)^2,
\end{equation}
to guarantee that $D(N,y(\overline{\delta}),x)<N/4$ for large $N$.
For the sake of simplicity, we take $J=2$ and the proposition is
proved as stated.
\endproof

\begin{Remark}
We remark here that condition \eqref{eq:conditionbeta} on $\beta$,
which we assumed to prove the proposition, was necessary to ensure the
existence of the function $\varphi$ satisfying all the properties we
needed, and in particular \eqref{eq:beta23}.
\end{Remark}

\begin{Corollary}
\label{cor:rho}
Let $A$ be a fixed positive constant and $\beta$ satisfy
\eqref{eq:conditionbeta}. Then $\vert\varrho(x)\vert<x^{-A}$ for all
sufficiently large $x\in\R$.
\end{Corollary}

\proof
It holds
\[
  \vert\varrho(x)\vert
  =
  \Biggl|
    \varrho_{\lfloor x \rfloor + 1}(x)
    \prod_{n > \lfloor x \rfloor + 1}\cos(\pi x a_n)
  \Biggr|
  < x^{-A}.
\qedhere
\]
\endproof

\begin{Theorem}
\label{th:general}
Let $C'>0$ satisfy \eqref{eq:C} and $\beta$ satisfy \eqref{eq:conditionbeta}. Then for all intervals $I\subseteq\mathbb{R}$ of length $\vert I\vert>\exp(-C'(\log N)^2)$ one has
\[
\mathbb{P}[X_N\in I]=\int_I g(x) \dx x+o(\vert I\vert),
\]
as $N\rightarrow\infty$, where
\[
g(x):=2\int_0^{\infty}\cos(2\pi u x)\prod_{n=1}^{\infty}\cos\left(\frac{2\pi u}{b_n}\right) \dx u
= 2\int_0^{\infty}\cos(2\pi u x) \varrho(2u) \dx u.
\]
\end{Theorem}

The proof follows along the same lines as Theorem 2.1 in
\cite{BettinMS2018b} and we omit the details for brevity.

\begin{Corollary}
\label{cor:general}
Let $\beta$ satisfy \eqref{eq:conditionbeta}. For all $\tau\in\R$ and $C'>0$ satisfying \eqref{eq:C}, we have
\[
\#\Set{(s_1,\dots,s_N)\in\Set{-1,+1}^N:\left\vert\tau-\sum_{n=1}^N\frac{s_n}{b_n}\right\vert < \delta} \sim 2^{N+1}g(\tau)\delta(1+o_{C',\tau}(1))
\]
as $N\to\infty$ and $\delta\to 0$, uniformly in $\delta\ge\exp(-C'(\log N)^2)$.
In particular, for large enough $N$, one has $\m_N(\tau)<\exp(-C'(\log N)^2)$.
\end{Corollary}

\begin{Remark}
\label{rk:beta}
We have imposed condition \eqref{eq:conditionbeta} for $\beta$ to keep the same range of validity for $x$ as in \cite{BettinMS2018b}. We remark that the hypotheses on $\beta$ could be relaxed at the price of restricting this range: for example, we could take
\[
\beta(N) = \frac{\log\log\log N}{\log\log N},
\]
and obtain the result of Proposition \ref{prop32} for $x\in
[N,\exp(\log^a N)]$, where $a\in (1,2)$ is a suitable
constant. In fact, this would weaken directly the estimates that we
have just found in Theorem~\ref{th:general} and Corollary
\ref{cor:general}, where $\exp(-C'(\log N)^2)$ would be replaced by
$\exp(-\log^a N)$.
\end{Remark}

\section{Products of \texorpdfstring{$k$}{k} primes}

We now leave the general case and concentrate on primes and products of $k$ distinct primes. Hence, we define
\[
\mathcal{P}_k := \Set{n\in\N | \text{$n$ is the product of $k$ distinct primes}};
\]
we will denote by $b_n^{(k)}$ the $n$-th element of the ordered set $\mathcal{P}_k$.
Let us recall the definition of $\mathcal{S}(N,\delta,x)$ in the case $a_n=1/b_n^{(k)}$:
\[
  \mathcal{S}(N,\delta,x):=\set{n\in\{1,\ldots,N\}:\Vert x/b_n^{(k)}\Vert\ge \delta}.
\]
We remark that, since we left the general case, we can now take $B(n) = b^{(k)}_n / n$, and denote it by $B_k(n)$. In 1900, Landau \cite{Landau1900} proved that
\[
\pi_k (t)
:= \vert\mathcal{P}_k \cap \Set{n\in\N | n\le t}\vert
= \frac{t}{\log t} \frac{(\log\log t)^{k-1}}{(k-1)!} + O\left(\frac{t(\log\log t)^{k-2}}{\log t}\right),
\]
which implies that
\begin{equation}
\label{eq:Bk}
B_k(n) \sim \log n \frac{(k-1)!}{(\log\log n)^{k-1}}.
\end{equation}
We can now start with a refinement of Proposition \ref{prop32}, where we extend the interval of validity for $x$ in the case $b_n = b_n^{(k)}$.

\begin{Proposition}
\label{stima_primi}
Let $A$ be a fixed positive constant, $k\in\N$ be fixed and $a_n = 1/ b_n^{(k)}$, where $b_n^{(k)}$ is the $n$-th element of the ordered set $\mathcal{P}_k$. Then $\vert\varrho_N(x)\vert<x^{-A}$ for all sufficiently large positive
integers $N$ and for all $x\in[U,\exp(f(N))]$, where $\log N = o(f(N))$ and
\[
f(N) = o\left(\left(\frac{N}{B_k^2(N)}\right)^{1/(2k+1)}\right),
\]
and $U>1$ is a constant depending on $f$.
\end{Proposition}

\proof
Let $x\in[N,\exp(f(N))]$. As in the proof of Proposition \ref{prop32}, we need to show that $D(N,y(\overline{\delta}),x)<N/4$, where $\overline{\delta}$ is chosen in the same way and $y(\overline{\delta}) = \overline{\delta} NB_k(N)$. Since now we are considering $x\ge N$, it is easy to see that for sufficiently large $N$ we have $y(\overline{\delta})\le x$.
We recall here that the prime omega function $\omega(n)$ is defined as
the number of different prime factors of $n$, and that
\[
\omega(n)\ll\frac{\log n}{\log \log n},
\]
as a consequence of the prime number theorem. In this case, we have
\begin{align*}
D(N,y(\overline{\delta}),x) &:=
\sum_{x-y(\overline{\delta})<m<x+y(\overline{\delta})}
\sum_{\substack{b_n^{(k)} \vert m \\ N/2\le n\le N}}1
\le \sum_{x-y(\overline{\delta})<m<x+y(\overline{\delta})}
\sum_{\substack{p_{j_1}\dots p_{j_k} \vert m \\ p_{j_i} \text{distinct primes}}}1 \\
&\le \sum_{x-y(\overline{\delta})<m<x+y(\overline{\delta})} \omega(m)^k
\le (2y(\overline{\delta})+1) \max_{m<x+y(\overline{\delta})}\omega(m)^k \\
&\ll (N\log x)^{1/2} B_k(N) \left(\frac{\log 2x}{\log\log 2x}\right)^k
\ll N^{1/2} B_k(N)\,(\log x)^{k+1/2},
\end{align*}
where we used the trivial bound for the prime omega function.
If we show that this quantity is $o(N),$ we are done. So we need
\[
\log x = o\left(\left(\frac{N}{B_k^2(N)}\right)^{1/(2k+1)}\right).
\]
Hence we can take any $f$ that satisfies
\[
f(N) = o\left(\left(\frac{N}{B_k^2(N)}\right)^{1/(2k+1)}\right),
\]
where we recall that $B_k$ satisfies \eqref{eq:Bk}.
The theorem is then proved for $x\in[N,\exp(f(N))]$. If $x<N$, it holds
\[
\vert\varrho_N(x)\vert \le \vert\varrho_{\lfloor x\rfloor}(x)\vert,
\]
hence the result we have just proved holds also whenever $x \le\exp\big(f(\lfloor x\rfloor)\big)$.
But there must exist $U >0$ such that this holds for any $x>U$, since $\log x = o(f(x))$.
\endproof

We are now ready to prove a more general version of Theorem 2.1 of \cite{BettinMS2018b} for the sequence $\big(b_n^{(k)}\big)_{n\in\N}$.
\begin{Theorem}
Let $f$ and $a_n$ be defined as in Proposition \ref{stima_primi}. Then for all intervals $I\subseteq\mathbb{R}$ of length $\vert I\vert>\exp(-f(N))$ one has
\[
\mathbb{P}[X_N\in I]=\int_I g(x) \dx x+o(\vert I\vert),
\]
as $N\rightarrow\infty$, where
\[
g(x):=2\int_0^{\infty}\cos(2\pi u x)\prod_{n=1}^{\infty}\cos\left(\frac{2\pi u}{b_n^{(k)}}\right) \dx u
= 2\int_0^{\infty}\cos(2\pi u x) \varrho(2u) \dx u.
\]
\end{Theorem}

\proof
The proof follows the one of Theorem 2.1 of \cite{BettinMS2018b}. Let $\varepsilon>0$ be fixed. We define
\begin{align*}
\xi= \xi_{N,-\varepsilon} := \exp(-(1-\varepsilon)f(N)),
\\
\xi_+= \xi_{N,+\varepsilon} := \exp(-(1+\varepsilon)f(N)), \\
&\xi_0 := \xi_{N,0} = \exp(-f(N)),
\end{align*} so that $\xi^{-1}<\xi_0^{-1}$ and Proposition \ref{stima_primi} holds for $x\in[N,\xi_0^{-1}]$. For an interval $I = [a,b]$ with $b-a>2\xi_0$, let us define $I^+ := [a-\xi,b+\xi]$ and $I^- := [a+\xi_+,b-\xi_+]$. Then one can construct two smooth functions $\Phi^\pm_{N,\varepsilon,I}(x):\mathbb{R}\to[0,1]$ (from now on, we will drop the subscripts when they are clear by the context) such that
\[
\begin{cases}
\supp \Phi^+ \subseteq I^+ & \\
\Phi^+(x)=1 &\text{for $x\in I$,} \\
\supp \Phi^- \subseteq I & \\
\Phi^-(x)=1 &\text{for $x\in I^-$,} \\
(\Phi^{\pm})^{(j)}(x)\ll_j \xi^{-j} &\text{for all $j\ge 0$.} \\
\end{cases}
\]
By the last equation, we know that the Fourier transforms of $\Phi^\pm$ satisfy
\begin{equation}
\label{phihat}
\widehat{\Phi^\pm}(x)\ll_B(1+\vert x\vert\xi)^{-B} \quad \text{for any $B>0$ and $x\in\R$.}
\end{equation}
Since
\[
\E[\Phi^{-}(X_N)]\le\mathbb{P}[X_N\in I]\le\E[\Phi^{+}(X_N)],
\]
we just need to show that
\[
\E[\Phi^{\pm}(X_N)]
= \int_{\mathbb{R}}\Phi^{\pm}(x)g(x) \dx x+o_\varepsilon(\vert I\vert).
\]
From now on, $\Phi$ will indicate either $\Phi^+$ or $\Phi^-$. By Lemma \ref{media_phi} we have
\[
\E[\Phi(X_N)]
= \frac12\int_{\mathbb{R}}\widehat{\Phi}(x/2)\varrho_N(x) \dx x
= I_1+I_2+I_3,
\]
where $I_1$, $I_2$ and $I_3$ are the integrals supported respectively in $\vert x\vert<N^{\varepsilon}$, $\vert x\vert\in[N^{\varepsilon}, \xi^{-(1+\varepsilon)}]$ and $\vert x\vert>\xi^{-(1+\varepsilon)}$. Note that $\xi^{-(1+\varepsilon)} = \exp((1-\varepsilon^2)f(N)) > \exp(\varepsilon \log N) = N^\varepsilon$, that $\xi^{-(1+\varepsilon)} = \xi_0^{-(1-\varepsilon^2)} < \xi_0^{-1}$, and that $\xi^{-(1+\varepsilon)}\cdot\xi = \xi^{-\varepsilon} = \xi_0^{-\varepsilon(1-\varepsilon)}\to +\infty$ as $N\to +\infty$.
By Lemma \ref{lem:3.2} and Corollary \ref{cor:rho}, we have
\begin{align*}
I_1 &= \frac{1}{2}\int_{-N^\varepsilon}^{N^{\varepsilon}}\widehat{\Phi}(x/2)\varrho_N(x) \dx x
= \frac{1}{2}\int_{-N^\varepsilon}^{N^{\varepsilon}}\widehat{\Phi}(x/2)\varrho(x) \dx x + O\left(\lVert\widehat{\Phi}\rVert_\infty N^{-1+3\varepsilon}\right) \\
&= \frac{1}{2}\int_\R\widehat{\Phi}(x/2)\varrho(x) \dx x
+ O_A\left(\lVert\widehat{\Phi}\rVert_\infty N^{-(A-1)\varepsilon}\right)
+ O\left(\lVert\widehat{\Phi}\rVert_\infty N^{-1+3\varepsilon}\right) \\
&= \int_\R\widehat{\Phi}(x)\varrho(2x) \dx x
+ O_\varepsilon\left(\lVert\Phi\rVert_1 N^{-1+3\varepsilon}\right),
\end{align*}
where to conclude we chose $A = A(\varepsilon)$ sufficiently large. For the second integral, we use Proposition \ref{stima_primi} and obtain
\begin{align*}
\vert I_2 \vert &\le \Vert\widehat{\Phi}\Vert_\infty \int_{N^\varepsilon}^{\xi^{-(1+\varepsilon)}} \vert\varrho_N(x)\vert \dx x
\le \Vert\Phi\Vert_1 \int_{N^\varepsilon}^{\xi^{-(1+\varepsilon)}} x^{-A} \dx x
\le \Vert\Phi\Vert_1 \int_{N^\varepsilon}^{+\infty} x^{-A} \dx x \\
&\ll_\varepsilon \Vert\Phi\Vert_1 N^{-A\varepsilon + \varepsilon}
\ll_\varepsilon \Vert\Phi\Vert_1 N^{-1},
\end{align*}
where, as before, to conclude we took $A=A(\varepsilon)$ sufficiently large.
For the last integral, we recall that trivially $\vert\varrho_N(x)\vert\le 1$; using the bound \eqref{phihat}, we obtain
\begin{align*}
\vert I_3 \vert &\le \int_{\vert x\vert > \xi^{-(1+\varepsilon)}}\vert\widehat{\Phi}(x/2)\vert\dx x
\ll_B \int_{\xi^{-(1+\varepsilon)}}^{+\infty}(1+ x \xi)^{-B} \dx x
= (B-1)(\xi^{-1}+\xi^{-(1+\varepsilon)})^{1-B} \\
&\ll_B \xi_0^{B-1}
= o_\varepsilon(\xi_0)
= o_\varepsilon(\vert I\vert),
\end{align*}
where to conclude we chose $B=B(\varepsilon)$ sufficiently large.
We can now put these results together: using Parseval's theorem and the fact that $\lVert\Phi\rVert_1 = O_\varepsilon(\vert I\vert)$, we get
\[
\E[\Phi(X_N)]
= \int_\R\widehat{\Phi}(x)\varrho(2x) \dx x
+ O_\varepsilon\left(\lVert\Phi\rVert_1 N^{-1+3\varepsilon}\right)
+ o_\varepsilon(\vert I\vert)
= \int_\R\Phi(x)g(x) \dx x + o_\varepsilon(\vert I\vert)
\]
and the theorem is then proved.
\endproof

\begin{Remark}
By Corollary \ref{cor:rho}, for any $n\in \N$ it holds
\[
\int_{-\infty}^{+\infty}\vert t^n \varrho(t)\vert\dx t <\infty,
\]
which implies by standard arguments (see e.g. \S 5 of \cite{Schmuland2003}) that the density $g$ is a smooth strictly positive function. Besides, by the same corollary, $g(x)\ll_D x^{-D}$ for any $D>0$.

\end{Remark}

\begin{Corollary}
For all $\tau\in\R$, we have
\[
\#\Set{(s_1,\dots,s_N)\in\Set{-1,1}^N:\left\vert\tau-\sum_{n=1}^N\frac{s_n}{b_n^{(k)}}\right\vert < \delta} \sim 2^{N+1}g(\tau)\delta(1+o_\tau(1))
\]
as $N\to\infty$ and $\delta\to 0$, uniformly in $\delta\ge\exp(-f(N))$, where $f$ is defined as in Proposition \ref{stima_primi}. In particular, for $N$ large enough, one has $\m_N(\tau)<\exp\big(-f(N)\big)$.
\end{Corollary}

\section{Addendum (by J. Benatar and A. Nishry): proof of Corollary \ref{BN}}

\begin{proof}
Let $c_m$ denote the $m$-th non-prime integer, so that $c_1 = 1$,
$c_2 = 4$, $c_3 = 6$, \dots
We first approximate $\tau$ with a restricted harmonic sum of the form
$\sum_{m \le M} s_m c_m$, where $M = M(N) = N - \pi(N)$.
Since $C_m := c_m / m \sim 1 $, we may apply Theorem \ref{genthm} to
obtain a sequence of signs $(s_n)_{n \le M} \in \left\{-1, 1 \right\}^M$
such that
\[
  -1 \leq \tau' := \sum_{m \le M} s_m c_m - \tau \leq 1.
\]
Moreover, taking $(p_n)_{n \in \N}$ to be the sequence of primes, we
have that $B(n) \sim \log n$ and hence we may apply
Theorem~\ref{primetheorem} to get a choice of signs
$(\sigma_n)_{n \le \pi(N)} \in \left\{-1,1 \right\}^{\pi(N)}$ such
that
\[
  \Bigl\vert \tau' - \sum_{n \le \pi(N)} \frac{\sigma_n}{p_n}
  \Bigr\vert
  \ll_{\tau, \varepsilon}
  \exp\Bigl( -N^{1/3-\varepsilon} \Bigr).
  \qedhere
\]
\end{proof}

\providecommand{\bysame}{\leavevmode\hbox to3em{\hrulefill}\thinspace}
\providecommand{\MR}{\relax\ifhmode\unskip\space\fi MR }
\providecommand{\MRhref}[2]{%
  \href{http://www.ams.org/mathscinet-getitem?mr=#1}{#2}
}
\providecommand{\href}[2]{#2}

\bigskip

\begin{tabular}{l}
Alessandro Gambini \\
Dipartimento di Matematica Guido Castelnuovo \\
Sapienza Universit\`a di Roma \\
Piazzale Aldo Moro, 5 \\
00185 Roma, Italia \\
email (AG): \texttt{alessandro.gambini@uniroma1.it} \\[10pt]
Remis Tonon, Alessandro Zaccagnini \\
Dipartimento di Scienze, Matematiche, Fisiche e Informatiche \\
Universit\`a di Parma \\
Parco Area delle Scienze, 53/a \\
43124 Parma, Italia \\
email (RT): \texttt{remis.tonon@unimore.it} \\
email (AZ): \texttt{alessandro.zaccagnini@unipr.it}
\end{tabular}

\newpage

\appendix

\section{Numerical data}

\begin{table}[h]

\begin{tabular}{|r|r|}
\hline
 $N$ & $\m_N(0) \cdot p_1\cdots p_N$ \\
\hline
 1 & 1 \cr
 2 & 1 \cr
 3 & 1 \cr
 4 & 23 \cr
 5 & 43 \cr
\hline
 6 & 251 \cr
 7 & 263 \cr
 8 & 21013 \cr
 9 & 1407079 \cr
10 & 4919311 \cr
\hline
11 & 818778281 \cr
12 & 2402234557 \cr
13 & 379757743297 \cr
14 & 3325743954311 \cr
15 & 54237719914087 \cr
\hline
16 & 903944329576111 \cr
17 & 46919460458733911 \cr
18 & 367421942920402841 \cr
19 & 17148430651130576323 \cr
20 & 1236225057834436760243 \cr
\hline
21 & 4190310920096832376289 \cr
22 & 535482916756698482410061 \cr
23 & 29119155169912957197310753 \cr
24 & 443284248908491516288671253 \cr
25 & 28438781483496930396689638231 \cr
\hline
26 & 10196503226925713726754541885481 \cr
27 & 137512198125317766267968137765087 \cr
28 & 5572821202475305606211985553786081 \cr
29 & 77833992457426020006787481021085581 \cr
30 & 24244850423688161715955346535954790877 \cr
\hline
31 & 2030349334778419995324119439659994086131 \cr
32 & 76860130392109667765387079377871685276909 \cr
33 & 5191970624445760882844533168270184721318637 \cr
34 & 329643209271348431895096550792159132283920307 \cr
35 & 19171590315567357340242017182966253037383120953 \cr
\hline
36 & 58192378490977430486851365332352874578233287403 \cr
37 & 837477642920747839191618216897250374978659503996169 \cr
38 & 130665466261033919414441892800025408642432364448372023 \cr
39 & 7541550169407232608689149525984967898398947805296216009 \cr
40 & 23868339955752715692132986729285170427530832996153507207 \cr
\hline
\end{tabular}
\bigskip

\caption{\label{Fig2}
The values, multiplied by $p_1\cdots p_N$, of the smallest signed harmonic sums with the first $N$ primes, with $N$ up to 40. See also \href{https://oeis.org/A332399}{OEIS A332399}.}

\end{table}

\begin{landscape}

\begin{table}

\begin{tabular}{|r|r|}
\hline
 $N$ & $\m_N(0) \cdot p_1\cdots p_N$ \\
\hline
41 & 3343165792500492306892396976512891068137770193474133826457 \cr
42 & 47233268931962642510303169511493601517566800154537867238057 \cr
43 & 93915329439868205746156163805290441755151986127947916375626793 \cr
44 & 50313439148416324581127610155641150127987318260569172331033593181 \cr
45 & 2035703788246113211455753014584246782664737720644793016891955087197 \cr
\hline
46 & 193768861589178044091624877468627581772116464350368833881209864412247 \cr
47 & 4664128549520402650533030541013467806288648880741654578068005845271177 \cr
48 & 252294099680710988063673862003152188841680135741161924018446904086039541 \cr
49 & 1641527055336324967995403445372629420483564255197731535006975381936073433 \cr
50 & 25436424505451332441928319474656471336874167655047366774702187882274894064063 \cr
\hline
51 & 1780024077761328763318128562703299120404666081323149178582620236480827415289259 \cr
52 & 115533643751466097619699345183033980786661230484621892531131629910924364040946261 \cr
53 & 34644520573176659229537081198934624126738529150336245449473941125320497104653817109 \cr
54 & 7369668963051661582966392617319633009625522375611294051784365401090471220946387592789 \cr
55 & 1999632582248468763357938742475072167566513418694128163881669512737786988287075374795317 \cr
\hline
56 & 151351981933638637742621357138936533979590998748883750430193460129876391573603481014628429 \cr
57 & 15302724902698188450027684974980553939987991074013402437579866232981371846926226684458406969 \cr
58 & 6269085432675155135477773589250562149563926327373176617473379555222137615792922214195964225281 \cr
59 & 429918790837116674905123858093668694474961832761345115366942177591943696826657060080682245858603 \cr
60 & 115809464188499233574522294110279752895686365776568444548440426304978721966632473743873345620708313 \cr
\hline
\end{tabular}

\bigskip

\caption{\label{Fig3}
The values, multiplied by $p_1\cdots p_N$, of the smallest signed harmonic sums with    the first $N$ primes, with $N$ between 41 and 60. See also \href{https://oeis.org/A332399}{OEIS A332399}.}

\end{table}

\end{landscape}

\begin{table}

\begin{tabular}{|r|r|}
\hline
 $N$ & $\Delta_N \cdot p_1\cdots p_N$ \cr
\hline
  1 & 1 \cr
  2 & 1 \cr
  3 & 1 \cr
  4 & 2 \cr
  5 & 22 \cr
\hline
  6 & 35 \cr
  7 & 263 \cr
  8 & 4675 \cr
  9 & 24871 \cr
 10 & 104006 \cr
\hline
 11 & 2356081 \cr
 12 & 6221080 \cr
 13 & 141769355 \cr
 14 & 6096082265 \cr
 15 & 6928889495 \cr
\hline
 16 & 367231143235 \cr
 17 & 1283811918935 \cr
 18 & 78312527055035 \cr
 19 & 5246939312687345 \cr
 20 & 372532691200801495 \cr
\hline
 21 & 8815359347599933286 \cr
 22 & 223849990729887044174 \cr
 23 & 6148176498383067879445 \cr
 24 & 179847837287937160817963 \cr
 25 & 663024394602752425373130 \cr
\hline
\end{tabular}

\bigskip

\caption{\label{Fig5}
The values, multiplied by $p_1\cdots p_N$, of the shortest distances $\Delta_N$ between different signed harmonic sums with the first $N$ primes, with $N$ up to~25.}

\end{table}


\begin{thebibliography}{10}

\bibitem{BettinMS2018a}
S.~Bettin, G.~Molteni, and C.~Sanna, \emph{Greedy approximations by signed
  harmonic sums and the {Thue}--{Morse} sequence}, Arxiv preprint 1805.00075,
  2018.

\bibitem{BettinMS2018b}
S.~Bettin, G.~Molteni, and C.~Sanna, \emph{Small values of signed harmonic sums}, C. R. Math. Acad. Sci.
  Paris \textbf{356} (2018), no.~11-12, 1062--1074.

\bibitem{BleicherE1975}
M.~N. Bleicher and P.~{Erd\H os}, \emph{The number of distinct subsums
  of {$\sum_1^N 1 / i$}}, Math. Comp. \textbf{29} (1975), 29--42.

\bibitem{BleicherE1976}
M.~N. Bleicher and P.~{Erd\H os}, \emph{Denominators of {Egyptian} fractions. {II}}, Illinois J. Math.
  \textbf{20} (1976), no.~4, 598--613.
  
\bibitem{CGTZ}
  M. Cafferata, A. Gambini, R. Tonon, and A. Zaccagnini, Sequence A332399 in The On-Line Encyclopedia of Integer Sequences (2020), published electronically at \url{https://oeis.org/A332399}.

\bibitem{Crandall2008}
R.~E. Crandall, \emph{Theory of {ROOF} walks}, Unpublished. Available at
  \url{http://www.reed.edu/physics/faculty/crandall/papers/ROOF11.pdf}, 2008.

\bibitem{Euler1737}
L.~Euler, \emph{Variae observationes circa series infinitas}, Commentarii
  academiae scientiarum imperialis Petropolitanae \textbf{9} (1737), 160--188.

\bibitem{Landau1900}
E.~Landau, \emph{Sur quelques probl\`emes relatifs \`a la distribution des nombres premiers},
  Bull. Soc. Math. France \textbf{28} (1900), 25--38.

\bibitem{Morrison1995}
K.~E. Morrison, \emph{Cosine products, {Fourier} transforms, and random
  sums}, Amer. Math. Monthly \textbf{102} (1995), no.~8, 716--724.

\bibitem{Morrison1998}
K.~E. Morrison, \emph{Random walks with decreasing steps}, Unpublished manuscript,
  California Polytechnic State University, 1998.

\bibitem{Ramanujan1915}
S.~Ramanujan, \emph{Highly composite numbers}, Proc. London Math. Soc.
  \textbf{14} (1915), 347--409.

\bibitem{Schmuland2003}
B.~Schmuland, \emph{Random harmonic series}, Amer. Math. Monthly
  \textbf{110} (2003), no.~5, 407--416.
  
\bibitem{Worley1976a}
R.~T. Worley, \emph{Signed sums of reciprocals. I}, J. Austral. Math. Soc. Ser. A \textbf{21} (1976), no.~4, 410--413.

\bibitem{Worley1976b}
R.~T. Worley, \emph{Signed sums of reciprocals. II}, J. Austral. Math. Soc. Ser. A \textbf{21} (1976), no.~4, 414--417.

\end{thebibliography}
\end{document}